\documentclass[11pt]{article}

\usepackage{amsmath}
\usepackage{graphicx}
\usepackage{algorithmic}
\usepackage{algorithm}
\usepackage{amssymb}
\usepackage[tight,footnotesize]{subfigure}
\usepackage[scaled]{helvet} 
\usepackage{fourier}
\usepackage{fullpage}

\usepackage[usenames]{color}
\def\PP{\mathbb{P}}
\def\E{\mathbb{E}}
\def\x{{\mathbf x}}

\def\P {{\mathbf{P}}}         
\def\Q {{\mathbf{Q}}}         
\def\CF{{\mathcal{S}}}
\def\psc{{\zeta}}
\def\pscI{{\zeta}}
\def\pscII{{\tilde{\zeta}}}
\def\CI{{C}}
\def\CII{{\tilde{C}}}
\def\CFI{{\mathcal{S}}}
\def\CFII{{\tilde{\mathcal{S}}}}
\newtheorem{thm}{Theorem}[section]

\newtheorem{lem}[thm]{Lemma}
\newtheorem{prop}[thm]{Proposition}
\newtheorem{defn}[thm]{Definition}

\def \remarks {\noindent {\bf Remarks.}\ \ }
\newcommand{\ceil}[1]{\left\lceil#1\right\rceil}
\newcommand{\norm}[1]{\left\Vert#1\right\Vert}
\newcommand{\abs}[1]{\left\vert#1\right\vert}
\newcommand\vect[1]{{\bf#1}}
\newcommand\matr[1]{{\bf#1}}

\newcommand\alphabf{{\boldsymbol{\alpha}}}
\newcommand{\argmin}{\operatornamewithlimits{argmin}}
\newcommand{\argmax}{\operatornamewithlimits{argmax}}

\newcommand{\Real}{\mathbb R}
\newcommand\RR[1]{\mathbb{R}^{#1}}

\newcommand{\dn}[1]{{}}
\newcommand{\rg}[1]{{}}
\DeclareMathOperator{\supp}{supp}
\DeclareMathOperator{\ext}{b-ext}

\DeclareMathOperator{\trace}{trace}
\DeclareMathOperator{\range}{range}

\begin{document}

\title{Near Oracle Performance and Block Analysis of Signal Space Greedy Methods} 

\author{R. Giryes and D. Needell}

\maketitle


\begin{abstract} Compressive sampling (CoSa) is a new methodology which
demonstrates that sparse signals can be recovered from a small number of linear
measurements.  Greedy algorithms like CoSaMP have been designed for this
recovery, and variants of these methods have been adapted to the case where
sparsity is with respect to some arbitrary dictionary rather than an orthonormal
basis.  In this work we present an analysis of the so-called Signal Space
CoSaMP method when the measurements are corrupted with mean-zero white Gaussian
noise.  We establish near-oracle performance for recovery of signals sparse in
some arbitrary dictionary.  In addition, we analyze the block variant of the method
for signals whose supports obey a block structure, extending the method into the 
model-based compressed sensing framework.  Numerical experiments
confirm that the block method significantly outperforms the standard method
in these settings. \end{abstract}

\section{Introduction}
\label{sec:intro}

We consider the compressive sensing problem which aims to recover a signal $\vect{x}\in\mathbb{R}^{d}$ from noisy measurements
\begin{eqnarray}
\label{eq:measurement}
\vect{y} = \matr{M}\vect{x}+ \vect{e},
\end{eqnarray}
where $\vect{M}\in \RR{m\times d}$ is a known linear operator and $\vect{e} \in \RR{d}$ is additive bounded noise, i.e. $\norm{\vect{e}}_2^2 \le \varepsilon^2$.  
A typical assumption in this context is that the signal $\vect{x}$ is \textit{sparse}.  There are several notions of sparsity, the simplest of which is that the signal itself has a small number of non-zero elements: $\norm{\vect{x}}_0 \le k$,
where $\norm{\vect{x}}_0 = |\supp(\vect{x})|$ denotes the $\ell_0$ quasi-norm.  We call such signals $k$-sparse.

A common approach to the compressive sensing problem utilizes the following optimization problem, deemed $\ell_1$-synthesis,
\begin{eqnarray}
\label{eq:l1_synthesis}
\hat{\vect{x}}_{\ell_1} = \argmin \norm{\vect{x}}_1 & s.t. & \norm{\vect{y} - \matr{M}\vect{x}}_2 \le \varepsilon.
\end{eqnarray}
One can guarantee accurate recovery using this approach when the measurement operator $\matr{M}$ satisfies the \textit{Restricted Isometry Property} (RIP)~\cite{Candes05Decoding}, which states that for some small enough constant $\delta_k < 1$,
$$
(1-\delta_{k})\norm{\vect{x}}^2 \leq \norm{\matr{M}\vect{x}}^2 \leq (1+\delta_k)\norm{\vect{x}}^2 \quad\text{for all $k$-sparse $\vect{x}$}.
$$
for some small enough constant $\delta_k < 1$.

It has been shown~\cite{Candes06Stable,Candes06Near, foucart10Sparse} that when the signal $\vect{x}$ is $k$-sparse and $\matr{M}$ satisfies the RIP with $\delta_{2k} < 0.4652$, the program~\eqref{eq:l1_synthesis} accurately recovers the signal,
\begin{eqnarray}
\label{eq:l1_rec_error}
\norm{\hat{\vect{x}}_{\ell_1} - \vect{x}}_2 \le C_{\ell_1}\varepsilon.
\end{eqnarray}

Another approach to solving the compressive sensing problem~\eqref{eq:measurement} is to use a greedy algorithm.  These methods typically identify elements of the support of the signal or estimate the signal iteration by iteration until some halting criterion is met.
Recently introduced methods that use this strategy are the CoSaMP \cite{Needell09CoSaMP},
  IHT \cite{Blumensath09Iterative}, and HTP \cite{Foucart11Hard} methods.  Greedy methods attempt to uncover the support of the signal iteratively, and then utilize a simple least-squares problem to estimate the entire signal.  
  
{\bfseries Sparsity in arbitrary dictionaries.}  Both the convex optimization and iterative methods provide rigorous recovery guarantees when the signal is sparse in some fixed orthonormal basis.  
However, this simple notion of sparsity limits the reality of compressive sensing applications, so we instead consider signals sparse in some dictionary $\matr{D} \in \RR{d \times n}$:
$$
\vect{x} = \matr{D}\alphabf \quad\text{for some}\quad \norm{\alphabf}_0 \le k.
$$
In this setting one can utilize the same $\ell_1$-synthesis program to obtain a candidate coefficient vector $\hat{\alphabf}_{\ell_1}$ and then estimate the signal $\vect{x}$ by $\hat{\vect{x}}_{\ell_1} = \matr{D}\hat{\alphabf}_{\ell_1}$.  Initial work on this problem shows that under stringent requirements on the dictionary $\matr{D}$, accurate recovery is possible~(see e.g.~\cite{Rauhut08Compressed,elad07Analysis}).  Recently, a sufficient and necessary dictionary-based null-space condition was derived by Chen et.al., which in particular shows that when the dictionary is full spark and highly coherent, the method fails~\cite{chen2013null}.

Alternatively, one can solve the $\ell_1$-analysis problem which minimizes coefficients in the analysis domain,
\begin{eqnarray}
\label{eq:l1_analysis}
\hat{\vect{x}}_{\ell_1} = \argmin \norm{\matr{D}^*\vect{x}}_1 & s.t. & \norm{\vect{y} - \matr{M}\vect{x}}_2 \le \varepsilon.
\end{eqnarray}
Here and throughout, the notation $\matr{D}^*$ denotes the adjoint of the matrix $\matr{D}$. 
In~\cite{Candes11Compressed}, the authors prove accurate recovery using this approach when the operators $\matr{M}$ and $\matr{D}$ satisfy the $\matr{D}$-RIP:
\begin{equation}\label{def:D_RIP}
(1-\delta_{k})\norm{\matr{D}\alphabf}^2 \leq \norm{\matr{M}\matr{D}\alphabf}^2 \leq (1+\delta_k)\norm{\matr{D}\alphabf}^2 \quad\text{for all $k$-sparse $\alphabf$}.
\end{equation}

Recently, the greedy approaches have also been adapted to the setting in which signals are sparse with respect to arbitrary dictionaries.  In particular, the \textit{Signal Space CoSaMP} variant \cite{Davenport13Signal} of the CoSaMP method~\cite{Needell09CoSaMP} is shown in Algorithm~\ref{alg:SSCoSaMP}.  Here and throughout, the subscript $T$ denotes the restriction to elements (or columns) indexed in $T$.  
$\CF_{\psc k}(\vect{z})$ denotes the operator which returns a support of size $\psc k$ that approximates the support of the best $k$-sparse representation of 
$\vect{z}$ in the dictionary $\matr{D}$, and $\matr{P}_T$ denotes the projection onto the range of $\matr{D}_T$.

\begin{algorithm}[ht]
\caption{Signal Space CoSaMP (SSCoSaMP)} \label{alg:SSCoSaMP}
\begin{algorithmic}[l]

\REQUIRE $k, \matr{M}, \matr{D}, \vect{y}, a$ where $\vect{y} = \matr{M}\vect{x}
+ \vect{e}$, $k$ is the sparsity of $\vect{x}$ under $\matr{D}$ and $\vect{e}$ is
the additive noise.  ${\CFII}_{\psc ak}$ and $\CFI_{\psc k}$ are two near optimal support selection schemes.

\ENSURE $\hat{\vect{x}}$: $k$-sparse approximation of
$\vect{x}$.

\STATE Initialize the support $T^0 =\emptyset$, the residual $\vect{y}_r^0 = \vect{y}$ and set $t = 0$.

\WHILE{halting criterion is not satisfied}

\STATE $t = t + 1$.

\STATE Find new support elements: $T_\Delta={\CFII}_{\psc ak}(\matr{M}^*\vect{y}^{t - 1}_r)$.

\STATE Update the support: $\tilde{T}^t = T^{t -1} \cup
T_\Delta$.

\STATE Compute the representation: $\vect{x}_p = \matr{D}(\matr{MD}_{\tilde{T}^{t}})^{\dag}\vect{y} = \matr{D}\left(\argmin_{\tilde\alphabf}\norm{\vect{y} -\matr{MD}\tilde\alphabf}_2^2 \text{ s.t. } \tilde\alphabf_{(\tilde{T}^t)^C}=0\right)$.

\STATE Shrink support: $T^t = \CFI_{\psc k}(\vect{x}_p)$.

\STATE Calculate new representation: ${\vect{x}}^t = \P_{T^t}\vect{x}_p$.

\STATE Update the residual:
$\vect{y}_r^t = \vect{y} - \matr{M}{\vect{x}}^t$.

\ENDWHILE
\STATE Form final solution $\hat{\vect{x}} = {\vect{x}}^t$.
\end{algorithmic}
\end{algorithm}

This method is analyzed in~\cite{Davenport13Signal}, under the assumption of the {$\matr{D}$-RIP}~\eqref{def:D_RIP} and the 
assumption that one has access to projections $\CF_k$ which satisfy
\begin{equation}\label{dnw_reqs}
\norm{\CF_k(\vect{z}) - \CF_k^{*}(\vect{z})}_2 \leq \min\left(c_1\norm{\CF_k^{*}(\vect{z})}_2, c_2\norm{\vect{z} - \CF_k^{*}(\vect{z})}_2 \right),
\end{equation}
where $\CF_k^{*}$ denotes the optimal projection:
\begin{eqnarray}\label{eq:optimal_sparse_projection}
\CF^*_{k}(\vect{z}) = \argmin_{\abs{T}\le k} \norm{\vect{z} - \P_T\vect{z}}_2^2.
\end{eqnarray} 
Under these requirements, the authors prove that the method accurately recovers the $k$-sparse signal 
as in~\eqref{eq:l1_rec_error}.   

Although the assumption on the approximate projections is also made for other 
methods~\cite{Blumensath11Sampling,Giryes13Greedy}, it is unknown whether such methods can be obtained.  
Recently, Giryes and Needell~\cite{Giryes14GreedySignal} relaxed these assumptions by introducing the 
notion of \textit{near-optimal projections}.
\begin{defn}
\label{def:C_optimal_proj_reg}
A pair of procedures ${\CFI}_{\pscI k}$ and ${\CFII}_{{\pscII} k}$ implies a pair of near-optimal projections $\P_{{\CFI}_{\pscI k}(\cdot)}$ and $\P_{{\CFII}_{{\pscII} k}(\cdot)}$ with constants $\CI_k$ and ${\CII}_k$ if for any $\vect{z} \in \Real^d$, $\abs{\CFI_{\pscI k}(\vect{z})} \le \pscI k$, with $\pscI \ge 1$, 
$\abs{{\CFII}_{{\pscII} k}(\vect{z})} \le {\pscII} k$, with ${\pscII} \ge 1$, and
\begin{eqnarray}
\label{eq:C_optimal_proj_reg}
\norm{\P_{{\CFII}_{{\pscII} k}(\vect{z})}\vect{z}}_2^2 \ge {\CII}_k\norm{\P_{\CF^*_k(\vect{z})}\vect{z}}_2^2\quad\text{as well as}\quad
 \norm{\vect{z}-\P_{\CFI_{\pscI k}(\vect{z})}\vect{z}}_2^2 \le  \CI_k\norm{\vect{z} - \P_{\CF^*_k(\vect{z})}\vect{z}}_2^2, 
\end{eqnarray}
where $\P_{\CF^*_k}$ denotes the optimal projection as in~\eqref{eq:optimal_sparse_projection}.
\end{defn}

We assume from this point onward that SSCoSaMP is run using such a pair of procedures.  
It has been proven in \cite{Giryes14GreedySignal}  that when the dictionary $\matr{D}$ is incoherent or satisfies the RIP, that many standard algorithms in compressive sensing give near-optimal projections satisfying~\eqref{eq:C_optimal_proj_reg}.  This improves upon previous results since even in this case, it is unknown whether any methods exist that satisfy the stricter requirements of~\eqref{dnw_reqs}.  In particular, they prove the following result:

\begin{thm}{\cite{Giryes14GreedySignal}}
\label{thm:general_bound_old}
Let $\matr{M}$ satisfy the $\matr{D}$-RIP~\eqref{def:D_RIP} with $\delta_{(3\pscI+1)k}$ ($\pscI\ge 1$).  Let $\CFI_{\pscI k}$ and $\CFII_{2\pscI k}$ be a pair of near optimal procedures (as in Definition~\ref{def:C_optimal_proj_reg}) with constants $\CI_k$ and ${\CII}_{2k}$.
Apply SSCoSaMP (with $a = 2$) and let ${\vect{x}}^t$ denote the approximation to the $k$-sparse signal $\vect{x}$ after $t$ iterations of Algorithm~\ref{alg:SSCoSaMP}.  If
$\delta_{(3\pscI+1) k}< \epsilon^2_{\CI_k,\CII_{2k},\gamma}$ and
 \begin{eqnarray}
\label{eq:C_k_tilda_C_2k_cond_old}
\left(1+\sqrt{\CI_k} \right)^2\left( 1-\frac{{\CII}_{2k}}{(1+\gamma)^2} \right)<1,
\end{eqnarray}
then after a constant number of iterations $t^*$ it holds that
\begin{eqnarray}
\label{eq:general_bound_old}
&& \hspace{-0.5in} \norm{{\vect{x}}^{t^*} -\vect{x}}_2 \le   \eta_0\norm{\vect{e}}_2,
\end{eqnarray}
where $\gamma$ is an arbitrary constant, and
$\eta_0$ is a constant depending on $\delta_{(3\pscI+1)k}$, $\CI_k$, ${\CII}_{2k}$ and $\gamma$.  The constant $\epsilon_{\CI_k,\CII_{2k},\gamma}$ is greater than zero if and only if \eqref{eq:C_k_tilda_C_2k_cond_old} holds.
\end{thm}

\subsection{Block sparsity}
Often, signals in practice have additional structure beyond simple sparsity.  One common model is that the support of the signal
is clustered together in one or more blocks.  This model accounts for the well-known observation that many significant signal
coefficients tend to reside close to one another, and is called \textit{block sparsity}.  Block sparsity appears in numerous compressed
sensing applications including DNA microarrays, magnetoencephalography, sensor networks and 
communication~\cite{stojnic2009reconstruction,eldar2009robust,baron2009distributed,wakin2005recovery,Blumensath09Sampling,baraniuk2010model}.  
The case when 
block sparsity is represented in an orthonormal basis, the literature in compressed sensing provides many algorithms (including a so-called
model-based CoSaMP method) and recovery guarantees~\cite{stojnic2009reconstruction,eldar2009robust,tropp2006algorithms,Blumensath09Sampling,baraniuk2010model}.  
In~\cite{baraniuk2010model}, Baraniuk et.al. define the block-sparse
restricted isometry property as follows.

\begin{defn}
The matrix $\matr{M}$ satisfies the block-sparse RIP (block-RIP) with constant $\delta$ if
$$
(1-\delta) \|\vect{x}\|_2^2 \leq \|\matr{M}\vect{x}\|_2^2 \leq (1+\delta)\|\vect{x}\|_2^2
$$
holds for all $k$-block-sparse signals $\vect{x}$.  That is, the above holds for all 
$\vect{x}\in S_{B,k}$ where 
$$ S_{B,k} := \{[\vect{x}_1^*, \ldots, \vect{x}_{d/B}^*]^* : \vect{x_i} \in \mathbb{R}^B, \vect{x_i} = 0 \text{ for } i\notin\Omega, |\Omega|=k \}.$$
\end{defn}

Under the assumption that the sampling matrix has small enough block-RIP constant $\delta$, it is shown that a modified CoSaMP
method offers robust block-sparse signal recovery~\cite{baraniuk2010model}\footnote{These results hold for more general models as well
as signal ensembles.  Here, we focus on the block-sparse model.}.
These results show that block-sparse signals can be recovered from
far fewer measurements than would have been required by a simple sparsity assumption alone.  Indeed, traditional results require
$m$ on the order of $Bk \log (d/(Bk))$ whereas the block-sparse model requires on the order of $Bk + k\log(d/k)$.  This is a 
significant reduction, especially when the block size $B$ is very large.  It is thus important to analyze and 
adapt signal space methods as well to this model.  In our setting, we consider block-sparse signal $\vect{x}$ with respect to the dictionary
$\matr{D}$ to be those in the set:
\begin{equation}\label{d-block-sparse}
S^{\matr{D}}_{B,k} := \{\vect{x} = \matr{D}\vect{\alphabf} : \vect{\alphabf} = [\alphabf_1, \ldots, \alphabf_{n/B}], \alphabf_i \in \mathbb{R}^B, \alphabf_i = 0 \text{ for } i\notin\Omega, |\Omega|=k \}.
\end{equation}
We call such signals $k$-block-sparse.  Analogously, we define the set of all such coefficient vectors by
\begin{equation}\label{d-block-sparse-coeff}
V^{\matr{D}}_{B,k} := \{\vect{\alphabf} : \vect{\alphabf} = [\alphabf_1, \ldots, \alphabf_{n/B}], \alphabf_i \in \mathbb{R}^B, \alphabf_i = 0 \text{ for } i\notin\Omega, |\Omega|=k \},
\end{equation}
and the set of all block-sparse  supports by
\begin{eqnarray}\label{d-block-sparse-supp}
W^{\matr{D}}_{B,k} = \left\{ T: \alphabf_T \in V_{B,k}^{\matr{D}} \right\}.
\end{eqnarray}
  Remark that for $B=1$ the $k$-block-sparse model coincides with the regular $k$-sparse framework.
We extend the block-RIP to the signal space setting.
\begin{defn}
\label{def:block_D_RIP}
The matrix $\matr{M}$ satisfies the block-sparse $\matr{D}$-RIP (block-\matr{D}-RIP) with constant $\delta_k$\footnote{We abuse notation and denote both the $\matr{D}$-RIP and the block-$\matr{D}$-RIP constants by $\delta_k$. The use will be clear from the context.}
\begin{equation}\label{block-d-rip}
(1-\delta_k) \|\vect{x}\|_2^2 \leq \|\matr{M}\vect{x}\|_2^2 \leq (1+\delta_k)\|\vect{x}\|_2^2
\end{equation}
holds for all $\vect{x} \in S^{\matr{D}}_{B,k}$.  
\end{defn}

  A straightforward adaption of Signal Space CoSaMP to block-sparse signals is described
by Algorithm~\ref{sscosamp-block}.  We add as input the block size $B$, and force the method to select blocks of coefficients at each
iteration.  The approximate projections used must now respect the block-sparsity structure; For this purpose we generalize Definition~\ref{def:C_optimal_proj_reg} for block-sparse signals.
\begin{defn}
\label{def:C_optimal_proj}
A pair of procedures ${\CFI}_{B,\pscI k}$ and ${\CFII}_{B,{\pscII} k}$ implies a pair of near-optimal projections $\P_{{\CFI}_{B,\pscI k}(\cdot)}$ and $\P_{{\CFII}_{B,{\pscII} k}(\cdot)}$ with constants $\CI_{k} = \CI_{B,k}$ and ${\CII}_{k} = {\CII}_{B,k}$\footnote{By abuse of notation we denote also the near optimality constants in the block case by $\CI_{k}$ and $\CII_{k}$. The value of $B$ will be clear from the context.}  if for any $\vect{z} \in \Real^d$, ${{\CFI}_{B,{\pscI} k}(\vect{z})} \in W^{\matr{D}}_{B,\pscI k}$, with $\pscI \ge 1$, 
${{\CFII}_{B,{\pscII} k}(\vect{z})}\in W^{\matr{D}}_{B,\pscII k}$, with ${\pscII} \ge 1$, and
\begin{eqnarray}
\label{eq:C_optimal_proj}
\norm{\P_{{\CFII}_{B,{\pscII} k}(\vect{z})}\vect{z}}_2^2 \ge {\CII}_{B,k}\norm{\P_{\CF^*_{B,k}(\vect{z})}\vect{z}}_2^2\quad\text{as well as}\quad
 \norm{\vect{z}-\P_{\CFI_{B,\pscI k}(\vect{z})}\vect{z}}_2^2 \le  \CI_{B,k}\norm{\vect{z} - \P_{\CF^*_{B,k}(\vect{z})}\vect{z}}_2^2, 
\end{eqnarray}
where ${\CF^*_k}$, the support selection procedure in the optimal projection, is given by
\begin{eqnarray}\label{eq:optimal_sparse_projection_block}
\CF^*_{B,k}(\vect{z}) = \supp\left(\argmin_{\vect{w} \in V^{\matr{D}}_{B,k}} \norm{\vect{z} - \matr{D}\vect{w}}_2^2\right).
\end{eqnarray} 
\end{defn}

\begin{algorithm}[ht]
\caption{Signal Space CoSaMP (SSCoSaMP) for block-sparse signals} \label{sscosamp-block}
\begin{algorithmic}[l]

\REQUIRE $k, B, \matr{M}, \matr{D}, \vect{y}, a$ where $\vect{y} = \matr{M}\vect{x}
+ \vect{e}$, $B$ is the block size, $k$ is the block-sparsity of $\vect{x}$ under $\matr{D}$ and $\vect{e}$ is
the additive noise.  ${\CFII}_{B,\psc ak}$ and ${\CFI}_{B,\psc k}$ are two near optimal support selection schemes that obey the block-sparsity constraint.

\ENSURE $\hat{\vect{x}}$: $k$-block-sparse approximation of
$\vect{x}$.

\STATE Initialize the support $T^0 =\emptyset$, the residual $\vect{y}_r^0 = \vect{y}$ and set $t = 0$.

\WHILE{halting criterion is not satisfied}

\STATE $t = t + 1$.

\STATE Find new support elements in blocks: $T_\Delta=\CFII_{B,\psc ak}(\matr{M}^*\vect{y}^{t - 1}_r)$.

\STATE Update the support: $\tilde{T}^t = T^{t -1} \cup
T_\Delta$.

\STATE Compute the representation: $\vect{x}_p = \matr{D}(\matr{MD}_{\tilde{T}^{t}})^{\dag}\vect{y} = \matr{D}\left(\argmin_{\tilde\alphabf}\norm{\vect{y} -\matr{MD}\tilde\alphabf}_2^2 \text{ s.t. } \tilde\alphabf_{(\tilde{T}^t)^C}=0\right)$.

\STATE Shrink support in blocks: $T^t = \CFI_{B,\psc k}(\vect{x}_p)$.

\STATE Calculate new representation: ${\vect{x}}^t = \P_{T^t}\vect{x}_p$.

\STATE Update the residual:
$\vect{y}_r^t = \vect{y} - \matr{M}{\vect{x}}^t$.

\ENDWHILE
\STATE Form final solution $\hat{\vect{x}} = {\vect{x}}^t$.
\end{algorithmic}
\end{algorithm}

\subsection{Our contribution}
In this work we extend the results of~\cite{Giryes14GreedySignal} to provide near-oracle recovery guarantees when the measurement noise $\vect{e}$ is mean-zero Gaussian noise.  We focus on the Signal Space CoSaMP method, but analogous results can be obtained for other methods.  Our main result is summarized by the following theorem.

\begin{thm}
\label{thm:general_bound}
Let $\vect{y} = \matr{M}\x +\vect{e}$, where $\matr{M}$ satisfies the $\matr{D}$-RIP~\eqref{def:D_RIP} with a constant $\delta_{(3\psc+1)k}$ ($\psc\ge 1$), $\x$ is a vector with a $k$-sparse representation under $\matr{D}$ and $\vect{e}$ is a white Gaussian noise with variance $\sigma^2$.  Let $\CFI_{\pscI k}$ and $\CFII_{2\pscI k}$ be a pair of near optimal procedures (as in Definition~\ref{def:C_optimal_proj_reg}) with constants $\CI_k$ and ${\CII}_{2k}$.
Apply SSCoSaMP (with $a = 2$) and let ${\vect{x}}^t$ denote the approximation to the $k$-sparse signal $\vect{x}$ after $t$ iterations of Algorithm~\ref{alg:SSCoSaMP}.  If $\delta_{(3\pscI+1) k}< \epsilon^2_{\CI_k,\CII_{2k},\gamma}$ and
 \begin{eqnarray}
\label{eq:C_k_tilda_C_2k_cond}
\left(1+\sqrt{\CI_k} \right)^2\left( 1-\frac{{\CII}_{2k}}{(1+\gamma)^2} \right)<1,
\end{eqnarray}
then after a constant number of iterations $t^*$ it holds with high probability that that
\begin{eqnarray}
\label{eq:general_bound}
 \hspace{-0.5in} \norm{{\vect{x}}^{t^*} -\vect{x}}_2^2 = O\left( k \log(n) \sigma^2 \right)  ,
\end{eqnarray}
where $\gamma$ is an arbitrary constant, and the constant $\epsilon_{\CI_k,\CII_{2k},\gamma}$ is greater than zero if and only if \eqref{eq:C_k_tilda_C_2k_cond} holds.
\end{thm}

\remarks

\noindent{\bfseries 1.}  This improves upon Theorem~\ref{thm:general_bound_old} in general, since $\norm{\vect{e}}_2$ is expected to be on the order of $\sqrt{n}\sigma$ when $\vect{e}$ is mean-zero Gaussian noise with variance $\sigma^2$.  These results align with those of standard compressive sensing when the dictionary $\matr{D}$ is the identity \cite{Candes07Dantzig, Bickel09Simultaneous,Giryes12RIP}.\\

\noindent{\bfseries 2.}  This bound is, up to a constant and a $\log(n)$ factor, the same as the one we get if we use an oracle that foreknows the true support of the original signal $\vect{x}$. The oracle estimator and its error will be defined and calculated hereafter. 
Note that the $\log$ factor is inevitable for any practical estimator that does not have access to oracle information \cite{Candes06Modern}. 

Our second contribution is to extend the analysis of Signal Space CoSaMP to the block-sparse setting.  We prove the following for
the recovery given by Algorithm~\ref{sscosamp-block}.

\begin{thm}
\label{thm:general_bound_block}
Let $\vect{y} = \matr{M}\x +\vect{e}$, where $\matr{M}$ satisfies the block-$\matr{D}$-RIP~\eqref{def:D_RIP} with a constant $\delta_{(3\psc+1)k}$ ($\psc\ge 1$), $\x$ is a vector with a block-$k$-sparse representation under $\matr{D}$ and $\vect{e}$ is a white Gaussian noise with variance $\sigma^2$.
Suppose that $\CFI_{B, \psc k}$ and $\CFII_{B, 2\psc k}$ are 
near optimal procedures (as in Definition~\ref{def:C_optimal_proj} with optimal projection~\eqref{eq:optimal_sparse_projection_block}) 
with constants $\CI_{k}$ and $\CII_{2k}$ respectively.
Apply SSCoSaMP (with $a = 2$) and let ${\vect{x}}^t$ denote the approximation to the block-$k$-sparse signal $\vect{x}$ after $t$ iterations of Algorithm~\ref{sscosamp-block}.  If $\delta_{(3\pscI+1) k}< \epsilon^2_{\CI_k,\CII_{2k},\gamma}$ and 
\eqref{eq:C_k_tilda_C_2k_cond} holds,
then after a constant number of iterations $t^*$ it holds with high probability that that
\begin{eqnarray}
\label{eq:general_bound_block}
 \hspace{-0.5in} \norm{{\vect{x}}^{t^*} -\vect{x}}_2^2 = O\left( k \log(n) \sigma^2 \right)  ,
\end{eqnarray}
where $\gamma$ is an arbitrary constant, and the constant $\epsilon_{\CI_k,\CII_{2k},\gamma}$ is greater than zero if and only if \eqref{eq:C_k_tilda_C_2k_cond} holds.
\end{thm}

\remarks\\

\noindent{\bfseries 1.} The condition of Theorem~\ref{thm:general_bound_block}  is more relaxed than the one of Theorem~\ref{thm:general_bound} as only $Bk +  k \log(n/k)$ measurements are needed for the block-$\matr{D}$-RIP to hold (See Theorem~3.3 in \cite{Blumensath09Sampling}), while for the regular $\matr{D}$-RIP, $Bk\log(n/(Bk))$ measurements are needed.\\

\noindent{\bfseries 2.} If $B=1$ then Theorem~\ref{thm:general_bound_block} coincides with Theorem~\ref{thm:general_bound}.

\subsection{Organization}
We establish some required notation and preliminary lemmas in Section~\ref{sec:notation}. In Section~\ref{sec:oracle} we present the oracle estimator in the signal domain and calculate its recovery error.
 In Section~\ref{sec:guarantees} we present our main results, which imply the near-oracle performance of Theorems~\ref{thm:general_bound} and \ref{thm:general_bound_block}.  
 Our proofs are included in Section~\ref{sec:proofs}.  We present numerical experiments for Algorithm~\ref{sscosamp-block} 
 in Section~\ref{exps} and conclude our work in Section~\ref{sec:discuss}.

\section{Notation and Consequences of Block-$\matr{D}$-RIP}
\label{sec:notation}

As usual, we let $\norm{\cdot}_2$ denote the Euclidean ($\ell_2$) norm of a vector, and $\norm{\cdot}$ the spectral ($\ell_2 \rightarrow \ell_2$) norm of a matrix.  We write the $d\times d$ identity matrix as $\matr{I}_d$. For an index set $T$, we denote by $\matr{D}_{T}$ the sub-matrix of $\matr{D}$ whose columns are indexed by $T$.  $\P_T = \matr{D}_T\matr{D}_T^\dag$ denotes the orthogonal projection onto $\range(\matr{D}_T)$ and 
$\Q_T = \matr{I}_d - \P_T$ the orthogonal projection onto its orthogonal
complement.

We next recall some elementary consequences of the block-$\matr{D}$-RIP, whose proofs are very similar to the ones of the $\matr{D}$-RIP and can be found in~\cite{Giryes13Greedy}.

\begin{lem}
\label{cor:MP_RIP_norm}
If $\;\matr{M}$ satisfies the block-$\matr{D}$-RIP with a constant $\delta_{k}$ then
\begin{eqnarray}
\label{eq:MP_RIP_norm}
\norm{\matr{M}\P_{T}}^2 \le 1+\delta_{k} \quad\text{and}\quad \norm{\P_{T}(\matr{I} - \matr{M}^*\matr{M})\P_{T}} \le \delta_k
\end{eqnarray}
for every $T \in W^{\matr{D}}_{B,k}$.
\end{lem}

\begin{lem}
\label{cor:D_RIP_norm_diff}
If $\;\matr{M}$ satisfies the block-$\matr{D}$-RIP~\eqref{block-d-rip} then
\begin{eqnarray}
\label{eq:D_RIP_norm_diff}
&& \norm{\P_{T_1}(\matr{I} - \matr{M}^*\matr{M})\P_{T_2}} \le \delta_k,  \nonumber
\end{eqnarray}
for any $T_1$ and $T_2$ with ${T_1}\in W^{\matr{D}}_{B,k_1}$, ${T_2} \in W^{\matr{D}}_{B,k_2}$, and $k_1+k_2\le k$.
\end{lem}

\begin{lem}[Approximate projections]
Let $\CFI_{B, \psc k}$ and $\CFII_{B, \psc k}$ be a pair of near-optimal procedures as in Definition~\ref{def:C_optimal_proj}.
For any vector $\vect{v}\in \Real^d$ that has a block-$k$-sparse representation, a support set ${T} \in W^{\matr{D}}_{B,k}$,
and any $\vect{z} \in \RR{d}$  we have that
\begin{eqnarray}
\label{eq:C_optimal_ineq}
&& \norm{\vect{z}-\P_{\CFI_{B, \psc k}(\vect{z})}\vect{z}}_2^2 \le \CI_k \norm{\vect{v}-\vect{z}}_2^2, \quad\text{and}\\
\label{eq:C_optimal_ineq_up}
&& \norm{\P_{\CFII_{B, \psc k}(\vect{z})}\vect{z}}_2^2 \ge {\CII}_k\norm{\P_{T}\vect{z}}_2^2.
\end{eqnarray}
\end{lem}

Finally, an elementary fact that we will also utilize.

\begin{prop}
\label{prop:norm2_ineq}
For any two given vectors $\vect{x}_1$, $\vect{x}_2$ and a constant $c>0$ it holds that
\begin{eqnarray}
\norm{\vect{x}_1+\vect{x}_2}_2^2 \le (1+c)\norm{\vect{x}_1}_2^2 + (1+\frac{1}{c})\norm{\vect{x}_2}_2^2.
\end{eqnarray}
\end{prop}

\section{The Oracle Estimator in the Signal Domain}
\label{sec:oracle}

Before we proceed to develop our main result for SSCoSaMP, we start by asking what is the error of an estimator that foreknows the support of the original signal $\vect{x}$.
 Let $T$ be the true support of $\vect{x}$, then the oracle estimator is  simply
\begin{eqnarray}
\label{eq:oracle_eq}
\hat{\vect{x}}_O = \matr{D}_T\left(\matr{M}\matr{D}_T \right)^\dag \vect{y},
\end{eqnarray}
i.e., the minimizer of
\begin{eqnarray}
\min_{\tilde{\vect{x}}} \norm{\vect{y} - \matr{M}\tilde{\vect{x}}}_2^2 & s.t. & \tilde{\vect{x}} = \matr{D}\tilde{\alphabf}, \tilde{\alphabf}_{T^C} = 0.
\end{eqnarray} 

The oracle's error is given by the following lemma 
\begin{lem}
Let $\matr{M}$ satisfy the $\matr{D}$-RIP \eqref{def:D_RIP} and $\vect{x}$ be a signal with a block-$k$-sparse representation $\alphabf$ under a dictionary $\matr{D}$.  Assume the measurements are corrupte with mean-zero Gaussian noise with variance $\sigma^2$.  The oracle estimator's error is 
\begin{eqnarray}
  \frac{Bk \sigma^2}{1+\delta_k}  \le \E\norm{\vect{x} - \hat{\vect{x}}_O }_2^2 \le \frac{Bk \sigma^2}{1-\delta_k}. 
\end{eqnarray} 
\end{lem}

{\em Proof:}
Since $\vect{x}$ is supported on $T$ we can write it as $\vect{x} = \matr{D} \alphabf =  \matr{D}_T \alphabf_T$.
Plugging \eqref{eq:measurement} in \eqref{eq:oracle_eq} we have 
\begin{eqnarray}
\hat{\vect{x}}_O = \matr{D}_T\left(\matr{M}\matr{D}_T \right)^\dag \left(\matr{M}\matr{D}_T\alphabf_T + \vect{e} \right) = \vect{x} + \matr{D}_T\left(\matr{M}\matr{D}_T \right)^\dag \vect{e}.
\end{eqnarray}
Thus, the oracle's error equals 
\begin{eqnarray}
\E\norm{\vect{x} - \hat{\vect{x}}_O }_2^2 = \E\norm{\matr{D}_T\left(\matr{M}\matr{D}_T \right)^\dag \vect{e} }_2^2.
\end{eqnarray}
Using the $\matr{D}$-RIP we have 
\begin{eqnarray}
\frac{1}{1+\delta_k}\E\norm{\matr{M}\matr{D}_T\left(\matr{M}\matr{D}_T \right)^\dag \vect{e} }_2^2
\le \E\norm{\vect{x} - \hat{\vect{x}}_O }_2^2 \le \frac{1}{1-\delta_k}\E\norm{\matr{M}\matr{D}_T\left(\matr{M}\matr{D}_T \right)^\dag \vect{e} }_2^2
\end{eqnarray}
The proof ends by noticing that $\matr{M}\matr{D}_T\left(\matr{M}\matr{D}_T \right)^\dag$ is a projection operator. Therefore, from the properties of white Gaussian noise we have
\begin{eqnarray}
\E\norm{\matr{M}\matr{D}_T\left(\matr{M}\matr{D}_T \right)^\dag \vect{e} }_2^2 = \trace\left(\matr{M}\matr{D}_T\left(\matr{M}\matr{D}_T \right)^\dag \right)\sigma^2 =  \trace\left(\left(\matr{M}\matr{D}_T \right)^\dag \matr{M}\matr{D}_T \right)\sigma^2 = Bk \sigma^2.
\end{eqnarray}
\hfill $\Box$ \bigskip

\section{Main Results}
\label{sec:guarantees}

Though the oracle's error is promising, it is unattainable as we do not have the support of the original signal. We turn to analyze the SSCoSaMP for block-sparse signals, which is a feasible algorithm for signal recovery.
We provide theoretical guarantees for its recovery performance when the measurement noise is Gaussian. We assume $a = 2$ in the algorithm, however, analogous results for other values can be obtained similarly.
We concentrate on the proof of SSCoSaMP for block-sparse signals.  In doing so, we also prove our main result for non-structured sparse signals, as that result is the special case when $B=1$. 

\subsection{Theorem Conditions} Before we present the proof of the main result,
we recall the conditions which guarantee the assumptions of
Theorem~\ref{thm:general_bound}.  The first requirement, that
$\delta_{(3\pscI+1) k}< \epsilon^2_{\CI_k,\CII_{2k},\gamma}$ for a constant
$\epsilon^2_{\CI_k,\CII_{2k},\gamma} >0$, holds for many families of random
matrices when $m \ge \frac{C}{\epsilon_k^2} \left( k \log(\frac{n}{k\epsilon_k})
+ Bk\log(\frac{1}{\epsilon_k}) \right)$ \cite{Candes06Near,
Rauhut08Compressed,Candes11Compressed, Blumensath09Sampling,
Mendelson08Uniform}. The more challenging assumption in the theorem is the
condition \eqref{eq:C_k_tilda_C_2k_cond}, which requires $\CI_k$ and
${\CII}_{2k}$ to be close to $1$. However, we do have an access to such
projection operators in many practical settings, and these are not supported by
the guarantees provided in previous results
\cite{Davenport13Signal,Blumensath11Sampling,Giryes13Greedy,Giryes13IHTconf}.
In fact, when the dictionary $\matr{D}$ is incoherent or satisfies the RIP
itself, then simple thresholding or standard compressive sensing algorithms can
be used for the projection.  See Sec. 4 of~\cite{Giryes14GreedySignal} for a
detailed discussion.

\subsection{SSCoSaMP Near-Oracle Guarantees}

As in~\cite{foucart10Sparse, Needell09CoSaMP}, our proof utilizes an iteration invariant which guarantees that each iteration exponentially reduces the recovery error, down to the noise floor.

\begin{thm} \label{thm:SSCoSaMP_iter_bound} 
Let $\vect{y} = \matr{M}\x +\vect{e}$, where $\matr{M}$ satisfies the block-$\matr{D}$-RIP~\eqref{def:D_RIP} with a constant $\delta_{(3\psc+1)k}$ ($\psc\ge 1$), $\x$ is a vector with a block-$k$-sparse representation under $\matr{D}$ and $\vect{e}$ is an additive noise vector. Suppose $\CFI_{B, \psc
k}$ and $\CFII_{B, 2\psc k}$ be a pair of near optimal projections as in
Definition~\ref{def:C_optimal_proj} with constants $\CI_k$
and $\CII_{2k}$ respectively. Then the estimate of SSCoSaMP, Algorithm~\ref{sscosamp-block}, at the $t$-th iteration satisfies
\begin{eqnarray}
\label{eq:SSCoSaMP_iter_bound} && \hspace{-0.3in} \norm{\vect{x}^t -\vect{x}}_2
\le \rho\norm{\vect{x} - \vect{x}^{t-1}}_2 +
\eta\norm{\P_{T_{\vect{e}}}\matr{M}^*\vect{e}}_2, \end{eqnarray} for constants
$\rho$ and $\eta$, and where \begin{eqnarray} \label{eq:T_e_def} T_{\vect{e}} =
\argmax_{\tilde{T}:\abs{\tilde{T}} \le 3\psc Bk}
\norm{\P_{\tilde{T}}\matr{M}^*\vect{e}}_2. \end{eqnarray} The iterates converge,
i.e. $\rho<1$, if $\delta_{(3\psc+1)k}< \epsilon^2_{\CI_k,\CII_{2k},\gamma}$,
for some positive constant $\epsilon^2_{\CI_k,\CII_{2k},\gamma}$, and
\eqref{eq:C_k_tilda_C_2k_cond} holds. \end{thm}

An immediate corollary of the above theorem yields the following.

\begin{thm} \label{cor:SSCoSaMP_bound} Let $\vect{y} = \matr{M}\x +\vect{e}$, where $\matr{M}$ satisfies the block-$\matr{D}$-RIP~\eqref{def:D_RIP} with a constant $\delta_{(3\psc+1)k}$ ($\psc\ge 1$), $\x$ is a vector with a block-$k$-sparse representation under $\matr{D}$ and $\vect{e}$ is a vector of additive noise. Suppose $\CFI_{B,\psc
k}$ and $\CFII_{B, 2\psc k}$ be a pair of near optimal projections as in
Definition~\ref{def:C_optimal_proj} with the optimal projection
\eqref{eq:optimal_sparse_projection_block} and with constants $\CI_k$
and $\CII_{2k}$ respectively.  Then after a
constant number of iterations $t^* =
\ceil{\frac{\log(\norm{\vect{x}}_2/\norm{\vect{e}}_2)}{\log(1/\rho)}}$ it holds
that \begin{eqnarray} && \hspace{-0.5in} \norm{\vect{x}^{t^*} -\vect{x}}_2 \le
\left(1 +  \frac{1-\rho^{t^*}}{1-
\rho}\right)\eta\norm{\P_{T_{\vect{e}}}\matr{M}^*\vect{e}}_2, \end{eqnarray}
where $\eta$ is a constant and $T_e$ is defined as in~\eqref{eq:T_e_def}.
\end{thm}

{\em Proof:}
By using \eqref{eq:SSCoSaMP_iter_bound} and recursion we have that after $t^*$ iterations
\begin{eqnarray}
&& \hspace{-0.3in} \norm{\vect{x}^{t^*} -\vect{x}}_2 \le \rho^{t^*}\norm{\vect{x} - \vect{x}^{0}}_2 +
  (1+\rho+\rho^2+\dots \rho^{t^*-1})\eta\norm{\P_{T_{\vect{e}}}\matr{M}^*\vect{e}}_2
  \\ \nonumber && \le\left(1 +  \frac{1-\rho^{t^*}}{1-\rho}\right)\eta\norm{\P_{T_{\vect{e}}}\matr{M}^*\vect{e}}_2,
\end{eqnarray}
where the last inequality is due to the equation of the geometric series, the choice of $t^*$, and the fact that $\vect{x}^0 =\vect{0}$.
\hfill $\Box$ \bigskip

To prove the near oracle bound we need the following lemma, whose proof is presented in Section~\ref{sec:proofs}. 
\begin{lem}\label{lem:oracle}
If $\matr{M}$ has the block-$\matr{D}$-RIP with a constant $\delta_{3\psc k}$ and $\vect{e}$ is zero-mean white Gaussian noise with variance $\sigma^2$ then with probability exceeding $ 1- \frac{2}{(3\psc k)!}n^{-\beta}$ we have
\begin{eqnarray}
\norm{\P_{T_{\vect{e}}}\matr{M}^*\vect{e}}_2 \le \sqrt{(1+\delta_{3\psc k})3\psc Bk}\left(1+\sqrt{2(1+\beta)\log(n)} \right)\sigma.
\end{eqnarray}
\end{lem}

This lemma together with Theorem~\ref{cor:SSCoSaMP_bound} provides the following near-oracle performance theorem.
\begin{thm}
\label{thm:SSCoSaMP_bound_near_oracle}
Assume the conditions of Theorem~\ref{thm:SSCoSaMP_iter_bound} and that $\vect{e}$ is a zero-mean white Gaussian noise with variance $\sigma^2$.  Then after a constant number of iterations $t^* = \ceil{\frac{\log(\norm{\vect{x}}_2/\norm{\vect{e}}_2)}{\log(1/\rho)}}$ it holds with probability exceeding $ 1- \frac{2}{(3\psc Bk)!}n^{-\beta}$ that
\begin{eqnarray}
&& \hspace{-0.5in} \norm{\vect{x}^{t^*} -\vect{x}}_2 \le
  \left(1 +  \frac{1-\rho^{t^*}}{1-\rho}\right)\eta \sqrt{(1+\delta_{3\psc k})3\psc Bk}\left(1+\sqrt{2(1+\beta)\log(n)} \right)\sigma.
\end{eqnarray}
\end{thm}

Note that Theorem~\ref{thm:SSCoSaMP_bound_near_oracle} implies our main result, Theorem~\ref{thm:general_bound}.  We have thus established that SSCoSaMP provides near-oracle performance when the noise is mean-zero Gaussian.

\section{Proofs}\label{sec:proofs}

\subsection{Proof of Lemma~\ref{lem:oracle}}
We rely on the proof technique of Lemma~3 in \cite{Donoho94idealdenoising}.
Without loss of generality, we prove for the case of $\sigma =1$. By simple scaling we get the above result for any value of $\sigma$. 
Using Lemma~\ref{cor:MP_RIP_norm} we have that for any $\vect{e}_1, \vect{e}_2 \in \RR{d}$ and any support $\tilde{T}$, $\abs{\tilde{T}} \le 3\psc Bk$,
\begin{eqnarray}
\norm{\P_{\tilde{T}}\matr{M}^*(\vect{e}_1- \vect{e}_2)}_2 \le 
\sqrt{1+\delta_{3\psc k}}\norm{\vect{e}_1- \vect{e}_2}_2.
\end{eqnarray}
Thus we can say that $\norm{\P_{\tilde{T}}\matr{M}^*\cdot}_2^2$ is a $\sqrt{1+\delta_{3\psc k}}$-Lipschitz functional.
Using trace and expectation properties we have
\begin{eqnarray}
\E \norm{\P_{\tilde{T}}\matr{M}^*\vect{e}}_2^2 = \E\left[\trace\left(\vect{e}^* \matr{M}\P_{\tilde{T}} \P_{\tilde{T}}\matr{M}^*\vect{e} \right)  \right]
=\trace\left(\matr{M}\P_{\tilde{T}} \P_{\tilde{T}}\matr{M}^* \E \left[\vect{e} \vect{e}^*\right] \right) = \trace\left(\matr{M}\P_{\tilde{T}} \P_{\tilde{T}}\matr{M}^*\right),
\end{eqnarray}
where the last equality is due to $\E\left[\vect{e} \vect{e}^*\right] = \matr{I}$. Note that $\trace(\matr{M}\P_{\tilde{T}} \P_{\tilde{T}}\matr{M}^*)$ equals the sum of the singular values of $\matr{M}\P_{\tilde{T}}$. Since $\P_{\tilde{T}}$ is a projection to a subspace of dimension $3\psc Bk$, there are at most $3\psc Bk$ non-zero singular values. 
By the block-$\matr{D}$-RIP, we thus have that 
\begin{eqnarray}
\label{eq:PTe_M_noise_mean_squared_error_ineqaulity}
\E \norm{\P_{\tilde{T}}\matr{M}^*\vect{e}}_2^2 \le (1+\delta_{3\psc k})3\psc Bk,
\end{eqnarray}
and from Jensen's inequality it follows that 
\begin{eqnarray}
\label{eq:PtildeT_M_noise_mean_error_ineqaulity}
\E \norm{\P_{\tilde{T}}\matr{M}^*\vect{e}}_2 \le \sqrt{(1+\delta_{3\psc k})3\psc Bk},
\end{eqnarray}
 
Using concentration of measure in Gauss space  \cite{Pisier86Probabilistic,Milman86Asymptotic} we have 
\begin{eqnarray}
\label{eq:P_PtildeT_M_concentration_of_measure_standard}
\PP\left( \left|\norm{\P_{\tilde{T}}\matr{M}^*(\vect{e})}_2 - \E\left[\norm{\P_{\tilde{T}}\matr{M}^*\vect{e}}_2 \right] \right| \ge t \right) \le 2\exp\left(-\frac{t^2}{2(1+\delta_{3\psc k})} \right).
\end{eqnarray}
Using \eqref{eq:PtildeT_M_noise_mean_error_ineqaulity} we have $\norm{\P_{\tilde{T}}\matr{M}^*\vect{e}}_2 - \sqrt{(1+\delta_{3\psc k})3\psc Bk}\le \norm{\P_{\tilde{T}}\matr{M}^*\vect{e}}_2 - \E\left[\norm{\P_{\tilde{T}}\matr{M}^*\vect{e}}_2 \right] $ and thus
\begin{eqnarray}
\label{eq:P_PtildeT_M_inequality_2}
\PP\left( \norm{\P_{\tilde{T}}\matr{M}^*\vect{e}}_2 - \sqrt{(1+\delta_{3\psc k})3\psc B k} \ge t \right) &\le & 
\PP\left( \norm{\P_{\tilde{T}}\matr{M}^*\vect{e}}_2 - \E\left[\norm{\P_{\tilde{T}}\matr{M}^*\vect{e}}_2 \right]  \ge t \right) \\ \nonumber  & \le & 
\PP\left( \left|\norm{\P_{\tilde{T}}\matr{M}^*\vect{e}}_2 - \E\left[\norm{\P_{\tilde{T}}\matr{M}^*\vect{e}}_2 \right] \right| \ge t \right).
\end{eqnarray}
Combining \eqref{eq:P_PtildeT_M_inequality_2} and \eqref{eq:P_PtildeT_M_concentration_of_measure_standard} yields
\begin{eqnarray}
\PP\left( \norm{\P_{\tilde{T}}\matr{M}^*\vect{e}}_2 \ge \sqrt{(1+\delta_{3\psc k})3\psc B k} + t \right) &\le & 2\exp\left(-\frac{t^2}{2(1+\delta_{3\psc k})} \right).
\end{eqnarray}
Selecting $t = \sqrt{(1+\delta_{3\psc k})3\psc B k}\sqrt{2(1+\beta)\log(n)}$ we have $e^{-\frac{t^2}{2(1+\delta_{3\psc k})} }= n^{ -3\psc B k (1+\beta)}$. Using a union bound we have
\begin{eqnarray}
&& \PP\left( \norm{\P_{T_{\vect{e}}}\matr{M}^*\vect{e}}_2 \ge \sqrt{(1+\delta_{3\psc k})3\psc Bk}\left(1+\sqrt{2(1+\beta)\log(n)} \right)\right) \\ \nonumber && \le 
\sum_{\tilde{T}:\abs{\tilde{T}} = 3\psc k}\PP\left( \norm{\P_{\tilde{T}}\matr{M}^*\vect{e}}_2 \ge \sqrt{(1+\delta_{3\psc k})3\psc Bk}\left(1+\sqrt{2(1+\beta)\log(n)} \right) \right)   \\ \nonumber && \le   2{n \choose 3\psc Bk}n^{-3\psc Bk (1+\beta)} \le \frac{2}{(3\psc Bk)!}n^{-\beta},
\end{eqnarray}
which completes the claim.
\hfill $\Box$ \bigskip

\subsection{Proof of Theorem~\ref{thm:SSCoSaMP_iter_bound}}
We turn now to prove the iteration invariant, Theorem~\ref{thm:SSCoSaMP_iter_bound}.  Instead of presenting the proof directly, we divide the proof into several lemmas.
The first lemma gives a bound for $\norm{\vect{x}_p -\vect{x} }_2$ as a function of $\norm{\P_{T_{\vect{e}}}\matr{M}^*\vect{e}}_2$ and $\norm{\Q_{\tilde{T}^t}(\vect{x}_p - \vect{x})}_2$.
\begin{lem}
\label{lem:SSCoSaMP_xp_bound}
If $\matr{M}$ has the $\matr{D}$-RIP with a constant $\delta_{3\psc k}$, then with the notation of Algorithm~\ref{alg:SSCoSaMP}, we have
\begin{eqnarray}
\label{eq:SSCoSaMP_xp_bound}
\norm{\vect{x}_p -\vect{x}}_2 \le \frac{1}{\sqrt{1-\delta_{(3\psc +1)k}^2}}\norm{\Q_{\tilde{T}^t}(\vect{x}_p - \vect{x})}_2 
+  \frac{1}{1-\delta_{(3\psc +1)k}}\norm{\P_{T_{\vect{e}}}\matr{M}^*\vect{e}}_2
\end{eqnarray}
\end{lem}
{\em Proof:}
Since $\vect{x}_p \triangleq \matr{D}\alphabf_p$ is the minimizer of $\norm{\vect{y} - \vect{M}\tilde{\vect{x}}}_2$ with the constraints $\tilde{\vect{x}} = \matr{D}\tilde{\alphabf}$ and $\tilde{\alphabf}_{(\tilde{T}^t)^C} =0$, then
\begin{eqnarray}
\langle \matr{M} \vect{x}_p - \vect{y}, \matr{M}\vect{v} \rangle =0
\end{eqnarray}
for any vector $\vect{v} = \matr{D}\tilde\alphabf$ such that $\tilde\alphabf_{(\tilde{T}^t)^C} =0$.
Substituting $\vect{y} = \matr{M}\vect{x} + \vect{e}$ with simple arithmetic gives
\begin{eqnarray}
\label{eq:xp_x_property}
\langle \vect{x}_p - \vect{x}, \matr{M}^*\matr{M}\vect{v} \rangle = \langle \vect{e}, \matr{M}\vect{v} \rangle
\end{eqnarray}
where $\vect{v} = \matr{D}\tilde\alphabf$ and $\tilde\alphabf_{(\tilde{T}^t)^C} =0$.
To bound $\norm{\P_{\tilde{T}^t}(\vect{x}_p - \vect{x})}_2^2$ we use \eqref{eq:xp_x_property} with $\vect{v} = \P_{\tilde{T}^t}(\vect{x}_p - \vect{x})$, which gives
\begin{align}\label{eq:Q_xp_x_norm}
 \norm{\P_{\tilde{T}^t}(\vect{x}_p - \vect{x})}_2^2 &= \langle \vect{x}_p - \vect{x}, \P_{\tilde{T}^t}(\vect{x}_p - \vect{x}) \rangle \\
&= \langle \vect{x}_p - \vect{x}, (\matr{I}_d - \matr{M}^*\matr{M})\P_{\tilde{T}^t}(\vect{x}_p - \vect{x}) \rangle + \langle \vect{e}, \matr{M}\P_{\tilde{T}^t}(\vect{x}_p - \vect{x}) \rangle\notag\\
&\le \norm{ \vect{x}_p - \vect{x}}_2 \norm{\P_{\tilde{T}^t \cup T} (\matr{I}_d - \matr{M}^*\matr{M})\P_{\tilde{T}^t}}_2 \norm{\P_{\tilde{T}^t}(\vect{x}_p - \vect{x})}_2  + \norm{ \P_{\tilde{T}^t}\matr{M}^*\vect{e}}_2\norm{ \P_{\tilde{T}^t} (\vect{x}_p - \vect{x})}_2\notag\\
&\le \delta_{(3\psc +1)k}\norm{ \vect{x}_p - \vect{x}}_2 \norm{\P_{\tilde{T}^t}(\vect{x}_p - \vect{x})}_2  + \norm{ \P_{\tilde{T}^t}\matr{M}^*\vect{e}}_2\norm{\P_{\tilde{T}^t}(\vect{x}_p - \vect{x})}_2,\notag
\end{align}
where the first inequality follows from the Cauchy-Schwartz inequality, the projection property that $\P_{\tilde{T}^t} = \P_{\tilde{T}^t}\P_{\tilde{T}^t}$ and the fact that $\vect{x}_p - \vect{x} = \P_{\tilde{T}^t \cup T}(\vect{x}_p - \vect{x})$.
The last inequality is due to the block-$\matr{D}$-RIP property, the facts that $\tilde{T}^t \in W^{\matr{D}}_{B, 3\psc k}$ and $T \in W^{\matr{D}}_{B, k}$ and Lemma~\ref{cor:D_RIP_norm_diff}.
After simplification of \eqref{eq:Q_xp_x_norm} by $\norm{\P_{\tilde{T}^t}(\vect{x}_p - \vect{x})}_2$ we have
\begin{eqnarray}
\nonumber \norm{\P_{\tilde{T}^t}(\vect{x}_p - \vect{x})}_2 \le \delta_{(3\psc +1)k}\norm{\vect{x}_p - \vect{x}}_2   + \norm{ \P_{\tilde{T}^t}\matr{M}^*\vect{e}}_2.
\end{eqnarray}
Utilizing the last inequality with the fact that 
$\norm{\vect{x}_p - \vect{x}}_2^2 = \norm{\Q_{\tilde{T}^t}(\vect{x}_p - \vect{x})}_2^2+ \norm{\P_{\tilde{T}^t}(\vect{x}_p - \vect{x})}_2^2$ gives
\begin{eqnarray}
&& \norm{\vect{x}_p - \vect{x}}_2^2  \le \norm{\Q_{\tilde{T}^t}(\vect{x}_p - \vect{x})}_2^2 + \left(\delta_{(3\psc +1)k}\norm{\vect{x}_p - \vect{x}}_2   +\norm{ \P_{\tilde{T}^t}\matr{M}^*\vect{e}}_2 \right)^2.
\end{eqnarray}
The last equation is a second order polynomial of $\norm{\vect{x}_p - \vect{x}}_2$. Thus its larger root is an upper bound for it and together with \eqref{eq:T_e_def} this gives the inequality in \eqref{eq:SSCoSaMP_xp_bound}. For more details look at the derivation of (13) in \cite{foucart10Sparse}.
\hfill $\Box$ \bigskip

The second lemma bounds $\norm{\vect{x}^{t} - \vect{x}}_2$ in terms of $\norm{\Q_{\tilde{T}^t}(\vect{x}_p - \vect{x})}_2$ and $\norm{\P_{T_{\vect{e}}}\matr{M}^*\vect{e}}_2$ using the first lemma.
\begin{lem}
\label{lem:SSCoSaMP_xt_bound1}
Given that $\CFI_{B, \psc k}$ is the first
 procedure as in Definition~\ref{def:C_optimal_proj} with  
constant $\CI_k$, if $\;\matr{M}$ has the $\matr{D}$-RIP with a constant $\delta_{(3\psc +1)k}$, then
\begin{eqnarray}
&& \hspace{-0.5in} \norm{\vect{x}^t -\vect{x}}_2 \le \rho_1\norm{\Q_{\tilde{T}^t}(\vect{x}_p - \vect{x})}_2+  \eta_1\norm{\P_{T_{\vect{e}}}\matr{M}^*\vect{e}}_2
\end{eqnarray}
\end{lem}
{\em Proof:}
We start with the following observation
\begin{eqnarray}
\label{eq:xt_x_diff_norm}
&& \hspace{-0.3in} \norm{\vect{x} -\x^t}_2 = \norm{\vect{x} - \vect{x}_p +  \vect{x}_p-\x^t}_2
 \le \norm{\vect{x}- \vect{x}_p}_2 + \norm{\x^t -  \vect{x}_p}_2,
\end{eqnarray}
where the last step is due to the triangle inequality. 
Using \eqref{eq:C_optimal_ineq} with the fact that $\x^t = \P_{\CFI_{B,\pscI k}( \vect{x}_p)} \vect{x}_p$ we have
\begin{eqnarray}
\label{eq:xt_x_diff_norm_coeff1}
\norm{\x^t -  \vect{x}_p}_2^2 \le \CI_k\norm{\x - \vect{x}_p}_2^2.
\end{eqnarray}
Plugging \eqref{eq:xt_x_diff_norm_coeff1} in \eqref{eq:xt_x_diff_norm} leads to 
\begin{align}
\label{eq:xt_x_diff_norm_all}
\norm{\x -\x^t}_2  &\le (1 + \sqrt{C_k})\norm{\vect{x} - \vect{x}_p}_2\le \frac{1 + \sqrt{C_k}}{\sqrt{1-\delta_{(3\psc +1)k}^2}}\norm{\Q_{\tilde{T}^t}( \vect{x}_p - \vect{x})}_2 
+  \frac{1 + \sqrt{C_k}}{1-\delta_{(3\psc +1)k}}\norm{\P_{T_{\vect{e}}}\matr{M}^*\vect{e}}_2
\end{align}
where for the last inequality we use Lemma~\ref{lem:SSCoSaMP_xp_bound}.
\hfill $\Box$ \bigskip

The last lemma bounds $\norm{\Q_{\tilde{T}^t}(\vect{x}_p - \vect{x})}_2$ with $\norm{\vect{x}^{t-1} - \vect{x}}_2$ and $\norm{\P_{T_{\vect{e}}}\matr{M}^*\vect{e}}_2$.
\begin{lem}
\label{lem:SSCoSaMP_Pxp_bound}
Given that $\CFII_{B,2\psc k}$ is a near optimal support selection scheme with a constant $\CII_{2k}$, if $\matr{M}$ has the block-$\matr{D}$-RIP with constants $\delta_{(3\psc +1)k}$ and $\delta_{2\psc k}$ then
\begin{eqnarray}
\label{eq:SSCoSaMP_Pxp_bound}
&& \hspace{-0.3in}\norm{\Q_{\tilde{T}^t}(\vect{x}_p - \vect{x})}_2 \le 
 \eta_2 \norm{\P_{T_{\vect{e}}}\matr{M}^*\vect{e}}_2
 + \rho_2\norm{\vect{x} - \vect{x}^{t-1}}_2.
\end{eqnarray}
\end{lem}
{\em Proof:}
Looking at the step of finding new support elements one can observe that $\P_{T_{\Delta}}$ is a near optimal projection operator for $\matr{M}^*\vect{y}^{t-1}_r = \matr{M}^*(\vect{y} - \matr{M}\vect{x}^{t-1})$.
Noticing that $T_{\Delta} \subseteq \tilde{T}^t$ and then using \eqref{eq:C_optimal_ineq_up} with $\P_{T^{t-1} \cup T}$ gives
\begin{eqnarray}
\label{eq:SSCoSaMP_PMy_Mx_ineq}
&& \hspace{-0.3in} \norm{\P_{\tilde{T}^t}\matr{M}^*(\vect{y} - \matr{M}\vect{x}^{t-1})}_2^2  \ge \norm{\P_{{T}_{\Delta}}\matr{M}^*(\vect{y} - \matr{M}\vect{x}^{t-1})}_2^2 \ge \CII_{2k}\norm{\P_{T^{t-1} \cup T}\matr{M}^*(\vect{y} - \matr{M}\vect{x}^{t-1})}_2^2.
\end{eqnarray}

We start by bounding the left hand side of \eqref{eq:SSCoSaMP_PMy_Mx_ineq} from above. Using Proposition~\ref{prop:norm2_ineq}
with $\gamma_1>0$ and $\alpha>0$ we have
\begin{eqnarray}
\label{eq:SSCoSaMP_PMy_Mx_ineq_lhs}
\norm{\P_{\tilde{T}^t}\matr{M}^*(\vect{y} - \matr{M}\vect{x}^{t-1})}_2^2 &\le &
(1+\frac{1}{\gamma_1})\norm{\P_{\tilde{T}^t}\matr{M}^*\vect{e}}_2^2
 + (1+\gamma_1)\norm{\P_{\tilde{T}^t}\matr{M}^*\matr{M}(\vect{x} - \vect{x}^{t-1})}_2^2
\\ \nonumber &\le & \frac{1+\gamma_1}{\gamma_1}\norm{\P_{\tilde{T}^t}\matr{M}^*\vect{e}}_2^2
 + (1+\alpha)(1+\gamma_1)\norm{\P_{\tilde{T}^t}(\vect{x} - \vect{x}^{t-1})}_2^2  \\ \nonumber &&
+(1+\frac{1}{\alpha})(1+\gamma_1)\norm{\P_{\tilde{T}^t}(\matr{I}_d-\matr{M}^*\matr{M})(\vect{x} - \vect{x}^{t-1})}_2^2
\\ \nonumber & \le &  \frac{1+\gamma_1}{\gamma_1}\norm{ \P_{\tilde{T}^t}\matr{M}^*\vect{e}}_2^2 - (1+\alpha)(1+\gamma_1)\norm{\Q_{\tilde{T}^t}(\vect{x} - \vect{x}^{t-1})}_2^2
 \\ \nonumber &&  +\left(1+\alpha +\delta_{(3\psc+1)k} + \frac{\delta_{(3\psc +1)k}}{\alpha}\right)(1+\gamma_1)\norm{\vect{x} - \vect{x}^{t-1}}_2^2,
\end{eqnarray}
where the last inequality is due to Lemma~\ref{cor:MP_RIP_norm} and \eqref{eq:D_RIP_norm_diff}.

We continue with bounding the right hand side of \eqref{eq:SSCoSaMP_PMy_Mx_ineq} from below.
For the first element we use Proposition~\ref{prop:norm2_ineq} with constants $\gamma_2>0$ and $\beta >0$, and \eqref{eq:D_RIP_norm_diff} to achieve
\begin{eqnarray}
\label{eq:SSCoSaMP_PMy_Mx_ineq_rhs1}
 \norm{\P_{{T^{t-1} \cup T}}\matr{M}^*(\vect{y} - \matr{M}\vect{x}^{t-1})}_2^2 & \ge & \frac{1}{1+\gamma_2}\norm{\P_{T^{t-1} \cup T}\matr{M}^*\matr{M}(\vect{x} - \vect{x}^t)}_2^2
-\frac{1}{\gamma_2}\norm{\P_{T^{t-1} \cup T}\matr{M}^*\vect{e}}_2^2
\\ \nonumber & \ge & \frac{1}{1+\beta}\frac{1}{1+\gamma_2}\norm{\vect{x} - \vect{x}^{t-1}}_2^2
-\frac{1}{\gamma_2}\norm{\P_{T^{t-1} \cup T}\matr{M}^*\vect{e}}_2^2
\\ \nonumber && -\frac{1}{\beta}\frac{1}{1+\gamma_2}\norm{\P_{T^{t-1} \cup T}(\matr{M}^*\matr{M} - \matr{I}_d)(\vect{x} - \vect{x}^{t-1})}_2^2
\\ \nonumber & \ge & (\frac{1}{1+\beta}-\frac{\delta_{(\psc+1)k}}{\beta})\frac{1}{1+\gamma_2}\norm{\vect{x} - \vect{x}^{t-1}}_2^2
-\frac{1}{\gamma_2}\norm{\P_{T^{t-1} \cup T}\matr{M}^*\vect{e}}_2^2.
\end{eqnarray}

By combining \eqref{eq:SSCoSaMP_PMy_Mx_ineq_lhs} and \eqref{eq:SSCoSaMP_PMy_Mx_ineq_rhs1} with
\eqref{eq:SSCoSaMP_PMy_Mx_ineq} and then using \eqref{eq:T_e_def} we have
\begin{eqnarray}
 (1+\alpha)(1+\gamma_1)\norm{\Q_{\tilde{T}^t}(\vect{x} - \vect{x}^{t-1})}_2^2  &
 \le &  \frac{1+\gamma_1}{\gamma_1}\norm{\P_{T_{\vect{e}}}\matr{M}^*\vect{e}}_2^2
 +\CII_{2k}\frac{1}{\gamma_2}\norm{\P_{T_{\vect{e}}}\matr{M}^*\vect{e}}_2^2
 \\ \nonumber &&
 + \left(1+\alpha +\delta_{(3\psc+1)k} + \frac{\delta_{(3\psc+1)k}}{\alpha}\right)(1+\gamma_1)\norm{\vect{x} - \vect{x}^{t-1}}_2^2
\\ \nonumber &&
- \CII_{2k}(\frac{1}{1+\beta}-\frac{\delta_{(1+\psc)k}}{\beta})\frac{1}{1+\gamma_2}\norm{\vect{x} - \vect{x}^{t-1}}_2^2.
\end{eqnarray}
Division of both sides by $(1+\alpha)(1+\gamma_1)$ yields
\begin{eqnarray}
 \norm{\Q_{\tilde{T}^t}(\vect{x} - \vect{x}^{t-1})}_2^2 & \le &   \bigg(\frac{1}{\gamma_1(1+\alpha)}
 +\frac{\CII_{2k}}{\gamma_2(1+\alpha)(1+\gamma_1)}
\bigg)\norm{\P_{T_{\vect{e}}}\matr{M}^*\vect{e}}_2^2
 \\ \nonumber && 
 + \bigg( 1 +\frac{\delta_{(3\psc+1)k}}{\alpha}
 - \frac{\CII_{2k}}{(1+\alpha)(1+\gamma_1)(1+\gamma_2)}(\frac{1}{1+\beta}-\frac{\delta_{(\psc+1)k}}{\beta}) \bigg)\norm{\vect{x} - \vect{x}^{t-1}}_2^2.
\end{eqnarray}
Substituting $\beta = \frac{\sqrt{\delta_{(\psc+1)k}}}{1- \sqrt{\delta_{(\psc+1)k}}}$ gives
\begin{eqnarray}
\norm{\Q_{\tilde{T}^t}(\vect{x} - \vect{x}^{t-1})}_2^2 & \le &   \bigg(\frac{1}{\gamma_1(1+\alpha)}
 +\frac{\CII_{2k}}{\gamma_2(1+\alpha)(1+\gamma_1)}
\bigg)\norm{\P_{T_{\vect{e}}}\matr{M}^*\vect{e}}_2^2
 \\ \nonumber &&
 + \bigg(1 +\frac{\delta_{(3\psc+1)k}}{\alpha}
 - \frac{\CII_{2k}}{(1+\alpha)(1+\gamma_1)(1+\gamma_2)}\left(1-\sqrt{\delta_{(\psc+1)k}}\right)^2\bigg)\norm{\vect{x} - \vect{x}^{t-1}}_2^2,
\end{eqnarray}
Using $  {\alpha = \frac{\sqrt{\delta_{(3\psc+1)k}}}{\sqrt{\frac{\CII_{2k}}{(1+\gamma_1)(1+\gamma_2)}}\left(1-\sqrt{\delta_{(\psc+1)k}}\right)- \sqrt{\delta_{(3\psc+1)k}}}}$ yields
\begin{eqnarray}
\norm{\Q_{\tilde{\Lambda}^t}(\vect{x} - \vect{x}^{t-1})}_2^2 & \le &    \bigg(\frac{1}{\gamma_1(1+\alpha)}
 +\frac{\CII_{2k}}{\gamma_2(1+\alpha)(1+\gamma_1)}
\bigg)\norm{\P_{T_{\vect{e}}}\matr{M}^*\vect{e}}_2^2
 \\ \nonumber &&
 + \left(-\bigg( \sqrt{\delta_{(3\psc+1)k}}-\sqrt{\frac{\CII_{2k}}{(1+\gamma_1)(1+\gamma_2)}}\left(1-\sqrt{\delta_{(\psc+1)k}}\right)
 \bigg)^2 +1\right)\norm{\vect{x} - \vect{x}^{t-1}}_2^2,
\end{eqnarray}
The values of $\gamma_1, \gamma_2$ give a tradeoff between the convergence rate and the size of the noise coefficient.
For smaller values we get better convergence rate but higher amplification of the noise.
We make no optimization on them and choose them to be $\gamma_1 = \gamma_2 = \gamma$ where $\gamma$ is an arbitrary number greater than $0$.
Thus we have
\begin{eqnarray}
\norm{\Q_{\tilde{T}^t}(\vect{x} - \vect{x}^{t-1})}_2^2 & \le &
   \bigg(\frac{1}{\gamma(1+\alpha)}
 +\frac{\CII_{2k}}{\gamma(1+\alpha)(1+\gamma)}
\bigg)\norm{\P_{T_{\vect{e}}}\matr{M}^*\vect{e}}_2^2
 \\ \nonumber && 
 + \bigg(-\bigg( \sqrt{\delta_{(3\psc+1)k}}-\frac{\sqrt{\CII_{2k}}}{1+\gamma}\left(1-\sqrt{\delta_{(\psc+1)k}}\right)
 \bigg)^2 +1\bigg)\norm{\vect{x} - \vect{x}^{t-1}}_2^2.
\end{eqnarray}
Using the triangle inequality and the fact that $\Q_{\tilde{T}^t}\vect{x}_p = \Q_{\tilde{T}^t}\vect{x}^{t-1} =0$ gives the desired result.

\hfill $\Box$ \bigskip

With the aid of the above three lemmas we turn to the proof of the iteration invariant, Theorem~\ref{thm:SSCoSaMP_iter_bound}.

{\em Proof of Theorem~\ref{thm:SSCoSaMP_iter_bound}:}
Substituting the inequality of Lemma~\ref{lem:SSCoSaMP_Pxp_bound} into the inequality of Lemma~\ref{lem:SSCoSaMP_xt_bound1} gives
\eqref{eq:SSCoSaMP_iter_bound} with $\rho = \rho_1\rho_2$ and $\eta = \eta_1 + \rho_1\eta_2$. The iterates converge if $\rho_1^2\rho_2^2 < 1$.
Since $\delta_{(\psc+1)k} \le \delta_{3\psc k} \le \delta_{(3\psc+1)k}$ this holds if
\begin{eqnarray}\label{interinequ}
&& \hspace{-0.3in} \frac{\left( 1+\sqrt{\CI_k} \right)^2 }{1-\delta_{(3\psc+1)k}^2}  
   \left(1-\left( \left(\frac{\sqrt{\CII_{2k}}}{1+\gamma} +1 \right)\sqrt{\delta_{(3\psc+1)k}} -\frac{\sqrt{\CII_{2k}}}{1+\gamma}
 \right)^2\right) < 1.
\end{eqnarray}
Since $\delta_{(3\psc+1)k} <1$, we have $\delta_{(3\psc+1)k}^2 < \delta_{(3\psc+1)k} $.  Using this fact and expanding~\eqref{interinequ} yields
the stricter condition 
\begin{eqnarray}
\label{eq:epsilon_quadratic_ineq}
&& \hspace{-0.6in} \left(1-\left(1+\sqrt{C_k}\right)^2 \left(\frac{\sqrt{\CII_{2k}}}{1+\gamma} +1 \right)^2\right){\delta_{(3\psc+1)k}} 
+
 2 \left(1+\sqrt{C_k}\right)^2\left(\frac{\sqrt{\CII_{2k}}}{1+\gamma} +1 \right)\frac{\sqrt{\CII_{2k}}}{1+\gamma}
 \sqrt{\delta_{(3\psc+1)k}} \\ \nonumber && 
 \hspace{3in} + 2\sqrt{\CI_k} + \CI_k   - \frac{\CII_{2k}(1+\sqrt{C_k})^2}{(1+\gamma)^2}< 0.
\end{eqnarray}
The above equation has a positive solution if and only if \eqref{eq:C_k_tilda_C_2k_cond} holds.
Denoting its positive solution by ${\epsilon_{\CI_k,\CII_{2k},\gamma}}$, we have that the expression holds when
$\delta_{(3\psc+1)k} \le \epsilon_{\CI_k,\CII_{2k},\gamma}^2$, which completes the proof.  Note that in the proof we have
\begin{eqnarray*}
 {\eta_1 =  \frac{1 + \sqrt{C_k}}{1-\delta_{(3\psc +1)k}}}, &&
{\eta_2^2 = \bigg(\frac{1+\delta_{3\psc k}}{\gamma(1+\alpha)}
 +\frac{(1+\delta_{(\psc +1)k})\CII_{2k}}{\gamma(1+\alpha)(1+\gamma)}\bigg)},\\
  {\rho_1^2 = \frac{\left(1 + \sqrt{C_k}\right)^2}{1-\delta_{(3\pscI+1)k}^2}}, &&
    {\rho_2^2 = 1-\bigg( \sqrt{\delta_{(3\psc+1)k}}-\frac{\sqrt{\CII_{2k}}}{1+\gamma}\left(1-\sqrt{\delta_{(\psc+1)k}}\right)
 \bigg)^2 }, \\ &&
   {\alpha = \frac{\sqrt{\delta_{(3\psc+1)k}}}{\sqrt{\frac{\CII_{2k}}{(1+\gamma_1)(1+\gamma_2)}}\left(1-\sqrt{\delta_{(\psc+1)k}}\right)- \sqrt{\delta_{(3\psc+1)k}}}}
   \end{eqnarray*}
     and $\gamma >0$ is an arbitrary constant.
\hfill $\Box$ \bigskip


\section{Numerical experiments with block-sparsity}\label{exps}

 \begin{algorithm}[t]
 \caption{$\epsilon$-Block Orthogonal Matching Pursuit} \label{alg:OMP_eps}
\begin{algorithmic}[l]

\REQUIRE $k, \matr{D}, \vect{z}$, where $\vect{z} = \vect{x}
+ \vect{e}$, $\vect{x} = \matr{D}\alphabf$, $\alphabf_0 \in V_{B,k}^{\matr{D}}$, where $\left\{ T_1, T_2, \dots, T_{n/B} \right\}$ is the set of the valid $1$-block-sparse supports,
and $\vect{e}$ is additive noise.

\ENSURE $\hat{\vect{x}}$: $k$-block-sparse approximation of
$\vect{x}$ supported on $\hat{T}$.

\STATE Initialize estimate $\hat{\x}^0 = \vect{0}$, residual $\vect{r}^0 = \vect{z}$, support $\hat{T}^0 =\check{T}^0 =\emptyset$
 and set $t = 0$.

\WHILE{$t \le k$}

\STATE $t = t + 1$.

\STATE New support element: $i^t = \argmax_{i : T_i \not\subset \in \check{T}^{t-1}}
\norm{\matr{D}^*_{T_i}\vect{r}^{t-1}}_2$.

\STATE Extend support: $\hat{T}^t = \hat{T}^{t-1} \cup T_{i^t}$.

\STATE Calculate a new estimate: $\hat{\vect{x}}^t = \matr{D}_{\hat{T}^t}\matr{D}_{\hat{T}^t}^\dag\vect{z}$.

\STATE Calculate a new residual: $\vect{r}^t = \vect{z} - \hat{\vect{x}}^t$.

\STATE Get maximal correlated column: $\hat{i}^t = \argmax_{i \in T_{i^t}}
\norm{\matr{D}^*_{i}\vect{r}^{t-1}}_2$

\STATE Support $\epsilon$-extension: $\check{T}^t = \check{T}^{t-1} \cup \ext_{\epsilon,2}\left(\left\{ \hat{i}^t \right\}\right)$.

\ENDWHILE

\STATE Set estimated support $\hat{T} = \check{T}^t$.

\STATE Form the final solution $\hat{\vect{x}} =  \matr{D}_{\hat{T}}\matr{D}_{\hat{T}}^\dag\vect{z}$.

\end{algorithmic}
\end{algorithm}
%

We perform similar experiments to the ones presented in  \cite{Davenport13Signal, Giryes14GreedySignal} for the overcomplete-DFT with redundancy factor $4$ and check the effect of the block sparsity. In \cite{Davenport13Signal, Giryes14GreedySignal} the signal coefficients were either all clustered together or all well separated.
Here we test the case of several separated clusters in the coefficients. 

We consider two setups: One with $k=2$ and $B=4$ and one with $k=4$ and $B=4$. 
We compare the performance of SSCoSaMP with OMP \cite{MallatZhang93}, $\epsilon$-OMP \cite{Giryes14GreedySignal,Giryes13OMP}, block-OMP (BOMP) \cite{Eldar10Block} and $\epsilon$-BOMP as the approximate projections.
We do not include other methods since a thorough comparison has been already performed in \cite{Davenport13Signal,Giryes14GreedySignal}, and the goal here is to check the effect of the block sparsity.

$\epsilon$-BOMP is an extension of BOMP, presented in Algorithm~\ref{alg:OMP_eps}. This algorithm uses the block-extension operator, which is a generalization to the extension operator in  \cite{Giryes13OMP}\footnote{In \cite{Giryes13OMP} it is referred to as $\epsilon$-closure but since closure bears a different meaning in mathematics we use a different name here.}. This operator extends the support to include also the indices of the blocks in the dictionary that contain at least one atom that is highly correlated with one of the atoms in the current support.

\begin{defn}[$\epsilon$-block-extension]
Let $0\le \epsilon < 1$, $\matr{D}$ be a fixed dictionary and  $\left\{ T_1, T_2, \dots, T_{n/B} \right\}$ be the set of the valid $1$-block-sparse supports.
The $\epsilon$-block-extension of a given support set $T$ is defined as
$$
\ext_{\epsilon,2}(T) = \left\{i \in T_q\; : \;\exists j\in T,\;\exists l\in T_q, ~ \frac{\abs{\langle\matr{d}_l, \matr{d}_j\rangle}^2}{\norm{\matr{d}_l}_2^2\norm{\matr{d}_j}_2^2}\ge  1-\epsilon^2\right\}.
$$
\end{defn}

\begin{figure}[htb]
\begin{minipage}[b]{.48\linewidth}
  \centering
  \centerline{\includegraphics[width=7.0cm]{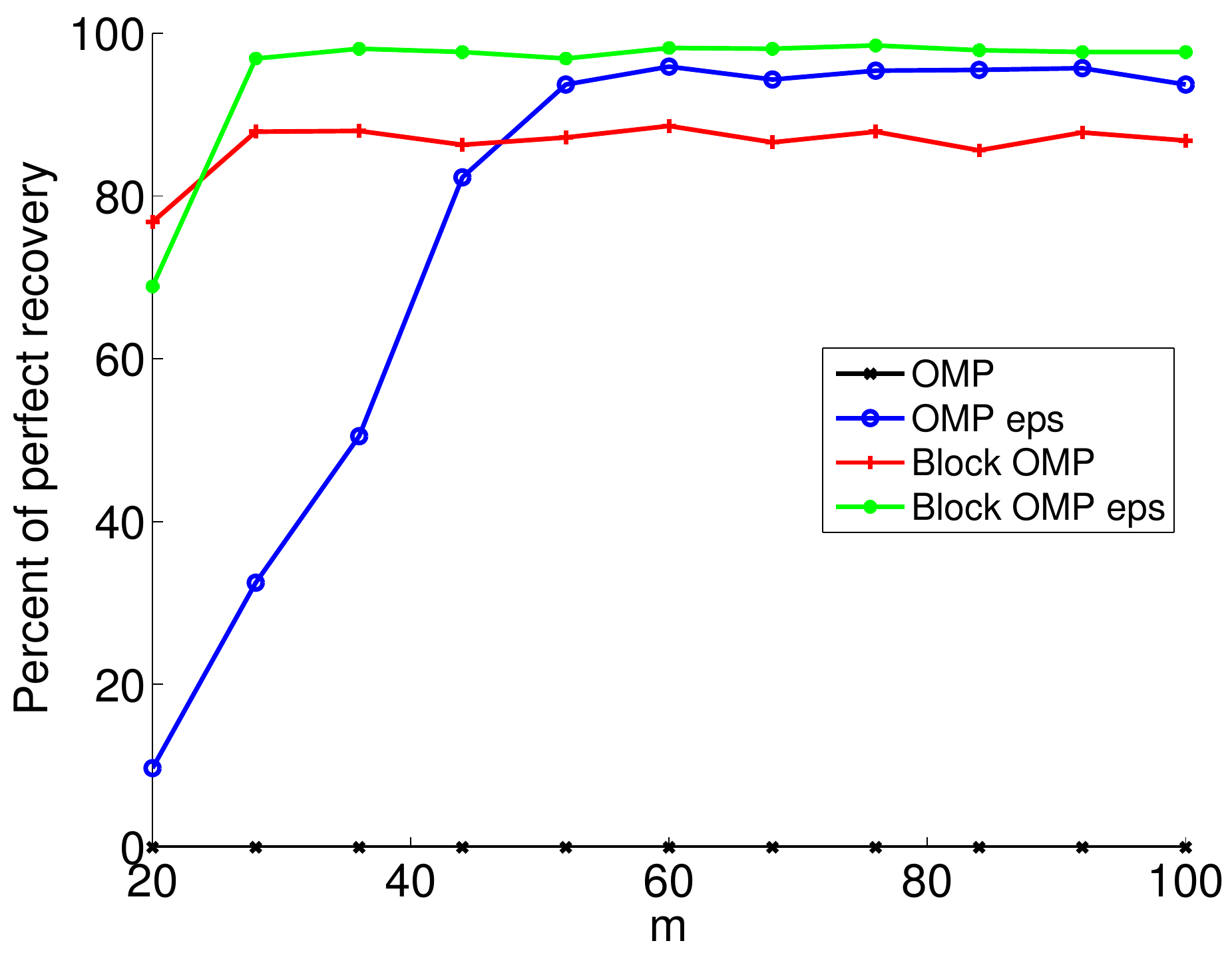}}
\end{minipage}
\hfill
\begin{minipage}[b]{.48\linewidth}
  \centering
  \centerline{\includegraphics[width=7.0cm]{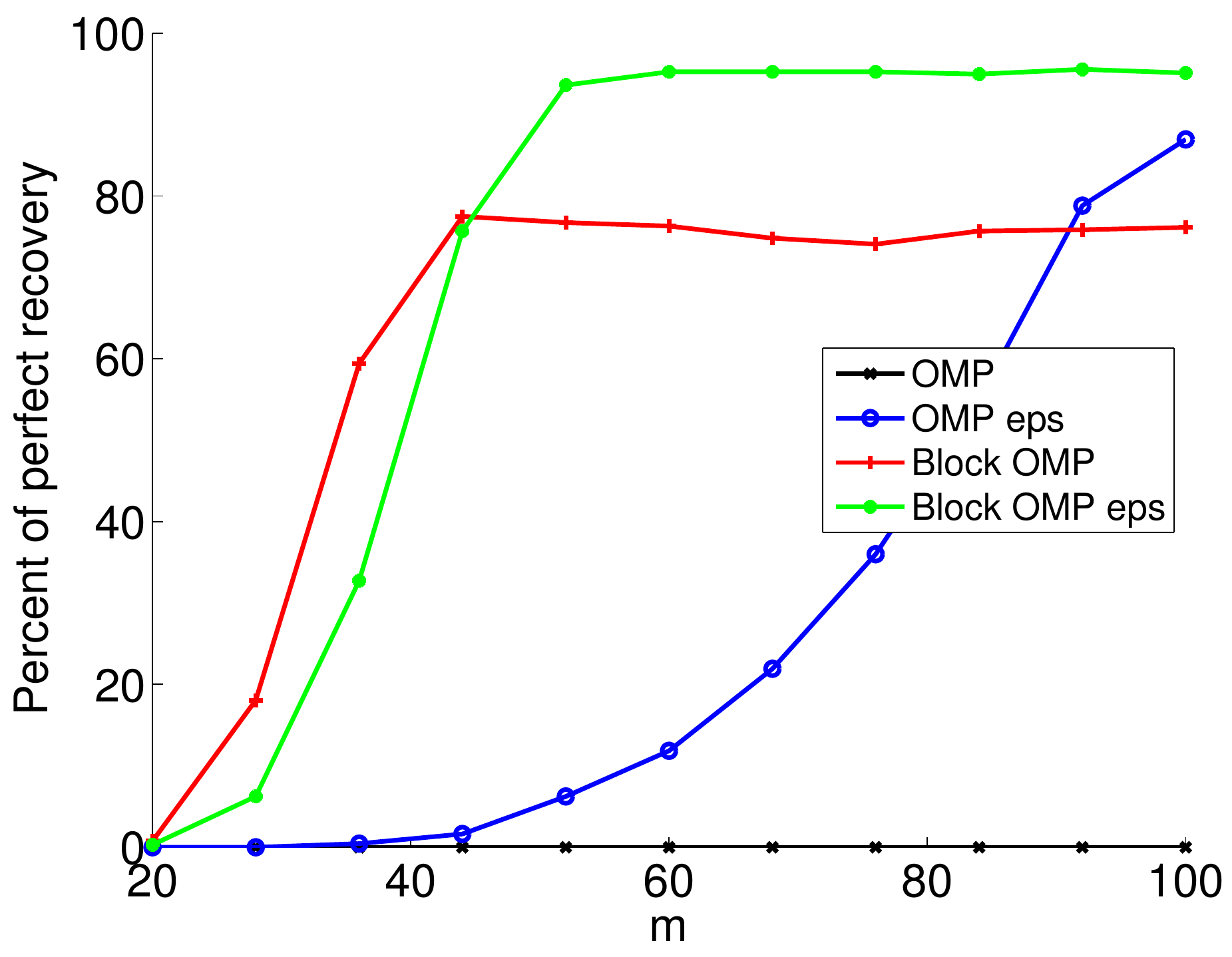}}
\end{minipage}
\caption{Recovery rate for 
SSCoSaMP (OMP), SSCoSaMP ($\epsilon$-OMP) with $\epsilon = \sqrt{0.1}$, SSCoSaMP (block-OMP) and SSCoSaMP ($\epsilon$-block-OMP) with $\epsilon = \sqrt{0.1}$ for a random $m\times 1024$ Gaussian matrix $\matr{M}$
and a $4$ times overcomplete DFT matrix $\matr{D}$.  The signal is $k$-block-sparse with block size $b=4$. On the left $k = 2$ and on the right $k=4$. 
}
\label{fig:recovery_rate}
\end{figure}

The recovery rate in the noiseless case appears in Figure~\ref{fig:recovery_rate}. It can be seen that using the block based algorithm provides better recovery results. Note that the major improvement is achieved for smaller values of $m$. This is due to the fact stated above that it is easier to satisfy the block-$\matr{D}$-RIP than the regular $\matr{D}$.
 
\begin{figure}[htb]
\begin{minipage}[b]{.48\linewidth}
  \centering
  \centerline{\includegraphics[width=7.0cm]{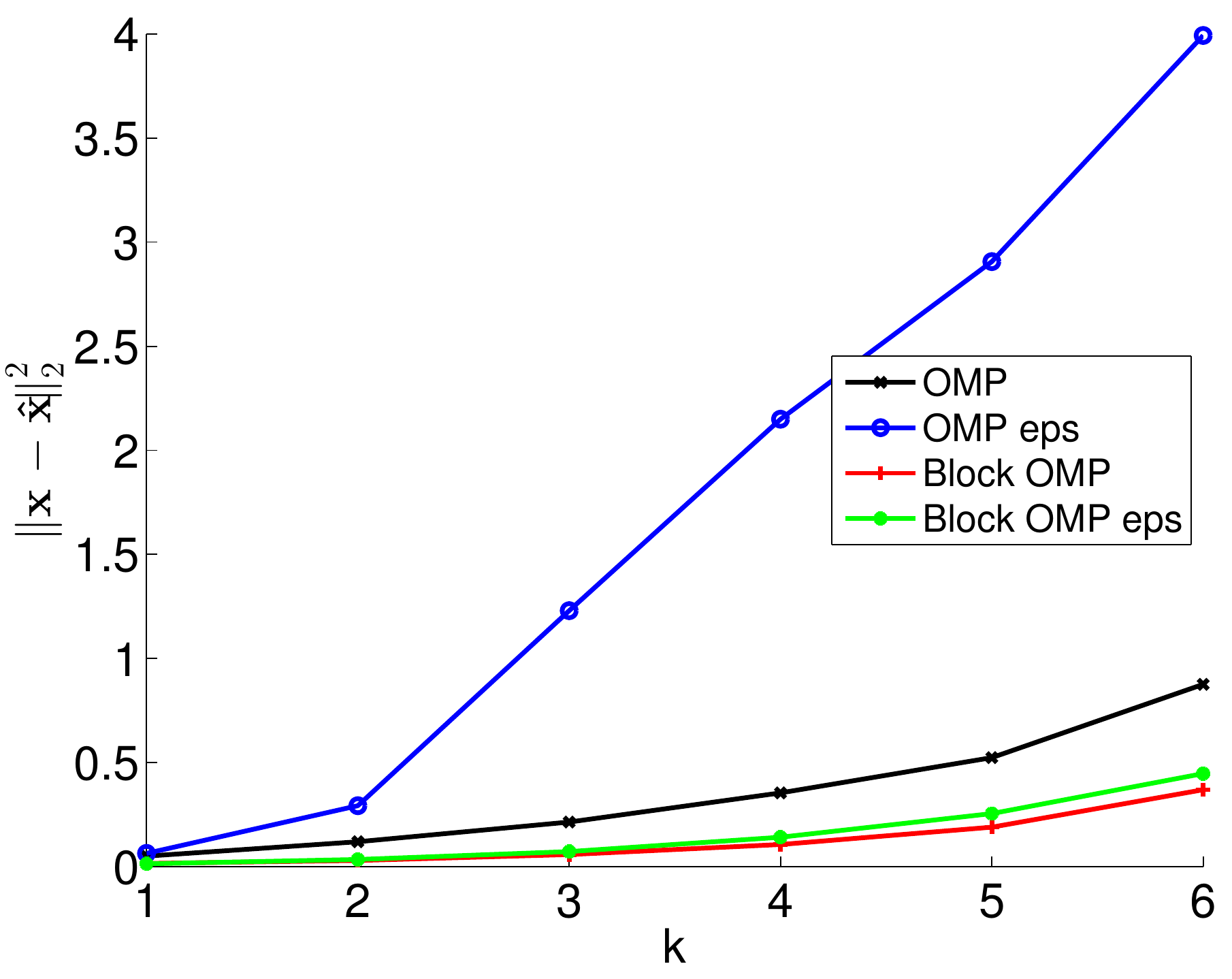}}
\end{minipage}
\hfill
\begin{minipage}[b]{.48\linewidth}
  \centering
  \centerline{\includegraphics[width=7.0cm]{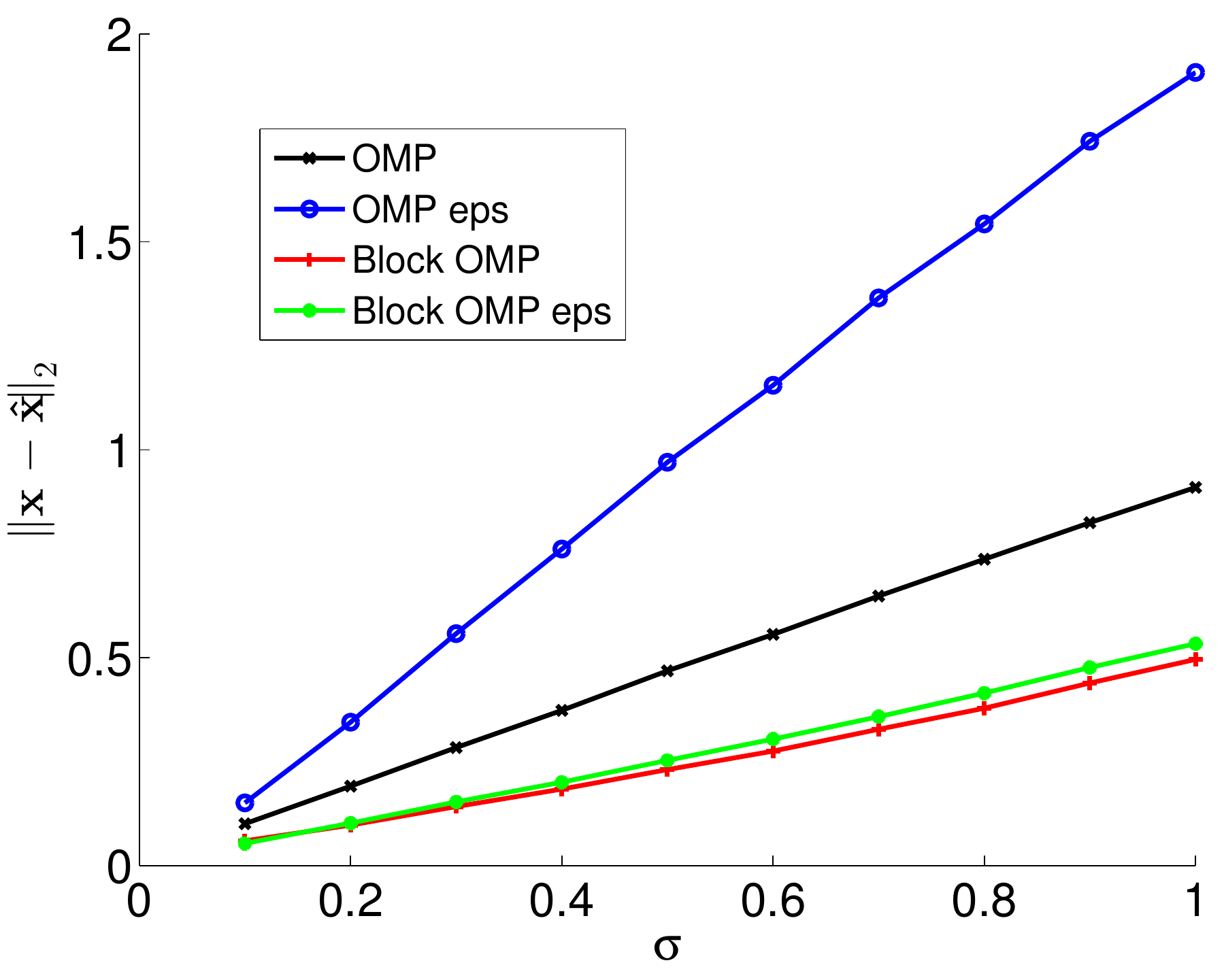}}
\end{minipage}
\caption{Error of 
SSCoSaMP (OMP), SSCoSaMP ($\epsilon$-OMP) with $\epsilon = \sqrt{0.1}$, SSCoSaMP (block-OMP) and SSCoSaMP ($\epsilon$-block-OMP) for a random $70\times 1024$ Gaussian matrix $\matr{M}$, a $4$ times overcomplete DFT matrix $\matr{D}$ and different values of $k$ with $\sigma = 0.4$ (left) and different values of $\sigma$ with $k= 2$. 
}
\label{fig:error_noisy}
\end{figure}

The reconstruction error in the noiseless case is presented in Figure~\ref{fig:error_noisy}. It can be seen that, as anticipated from the theory, the $\ell_2$-norm of the error scales linearly with $\sigma$ and that the squared error scales nearly linearly with $k$. The reason that the behavior for $k$ is not exactly linear as is for $\sigma$ is due to the fact that the $\matr{D}$-RIP constant changes when $k$ increases.  Notice also that the larger $k$ is, the harder it is to satisfy the $\matr{D}$-RIP conditions.  As the conditions for SSCoSaMP without using block-sparsity are stricter we have that the graphs in this case rise faster.

\section{Conclusion}
\label{sec:discuss}

 The Signal Space CoSaMP method was studied in the case of arbitrary noise~\cite{davenport2012compressive,Davenport13Signal,Giryes13Greedy} 
 under the assumptions of the $\matr{D}$-RIP and approximate projections.  As in~\cite{Giryes13Greedy}, the assumptions in this work on 
 the approximate projections allow for standard compressed sensing algorithms to be used when the dictionary $\matr{D}$ satisfies properties 
 like the RIP or incoherence.  In this correspondence, we have presented performance guarantees for this method in the presence of white 
 Gaussian noise, which are comparable to those obtained from an oracle which provides the support of the signal.  Our bounds are also of 
 the same order as those for standard greedy algorithms like IHT and CoSaMP~\cite{Giryes12RIP}, but ours hold also for signals sparse 
 with respect to an arbitrary dictionary.   
 
 In addition, we present a block variant of the Signal Space CoSaMP method designed for signals sparse in an arbitrary dictionary
 $\matr{D}$ and whose sparsity pattern obeys a block structure.  The analysis demonstrates that far fewer measurements are
 required when this model is utilized when compared to standard methods.  Experiments show that using traditional compressed sensing
 algorithms as the approximate projections in the Signal Space CoSaMP 
 method very often fails completely for block-sparse signals.  The block variant proposed in this work however, offers recovery
 in these settings. 


\section*{Acknowledgment}
R. Giryes is grateful to the Azrieli Foundation for the award of
an Azrieli Fellowship.  D. Needell was partially supported by the
Simons Foundation Collaboration grant $\#274305$, the Alfred P. Sloan fellowship, 
and NSF Career grant $\#1348721$.  In addition, the authors thank the reviewers
of the manuscript for their suggestions which greatly improved the paper.

\bibliographystyle{elsarticle-num}
\bibliography{SSCoSaMP,IEEEabrv}

\end{document}